%% file: diophantine_merged_siam.tex
\documentclass[]{siamart0516}
\usepackage[utf8]{inputenc}
\usepackage{pstricks,pst-node,pst-text,pst-3d,amssymb,amsmath,amsfonts,amscd,bm,epsfig,relsize,subfigure,fixltx2e,graphicx,algorithm}
\usepackage[noend]{algpseudocode}
\makeatletter
\usepackage[]{ntheorem}
\renewtheoremstyle{plain}%
{\item[\hskip\labelsep \theorem@headerfont ##1\ ##2\theorem@separator]}%
{\item[\hskip\labelsep \theorem@headerfont ##1\ ##2, ##3\theorem@separator]}
\makeatother
\newtheorem{exmp}{Example}

\usepackage{setspace}
\usepackage{array}
\usepackage{watermark}
\newcommand {\tr} {\rm\scriptscriptstyle T}

\newcommand{\bea}[1]{\begin{eqnarray} #1 \end{eqnarray}}
\usepackage{threeparttable, caption}
\newcommand {\eqrefn}{Eq.~\eqref}
\newcommand{\rank}[1]{\operatorname{rank}\Big(#1\Big)}
 \newcommand{\bbf}{\mathbf{b}}

\newcommand{\Abf}{\mathbf{A}}

\newcommand{\dbl}{\boldsymbol{d}}  
 \newcommand{\hbl}{\boldsymbol{h}}

  \newcommand{\xbl}{\boldsymbol{x}}
\newcommand{\ybl}{\boldsymbol{y}} \newcommand{\zbl}{\boldsymbol{z}}

\newcommand{\Abl}{\boldsymbol{A}}  
 \newcommand{\Ebl}{\boldsymbol{E}} 
\newcommand{\Gbl}{\boldsymbol{G}} \newcommand{\Hbl}{\boldsymbol{H}}

 \newcommand{\Qbl}{\boldsymbol{Q}} \newcommand{\Rbl}{\boldsymbol{R}}
  
  \newcommand{\Xbl}{\boldsymbol{X}}
\newcommand{\Ybl}{\boldsymbol{Y}} 
\newcommand{\bMatsmall}[1]{\left[\begin{smallmatrix}#1\end{smallmatrix}\right]}
\newcommand{\bMat}[1]{\begin{bmatrix}#1\end{bmatrix}}

\newcommand{\zeros}{\textbf{0}}

\def\BState{\State\hskip-\ALG@thistlm}

\newcommand{\Rcal}{\mathcal{R}}
 
 \newcommand{\Zcal}{\mathcal{Z}}
 \newcommand{\Scal}{\mathcal{S}}
{\theoremstyle{break}\newtheorem{thm}{Theorem}\theoremheaderfont{\normalfont\bfseries}}
{\theoremstyle{empty}\theorembodyfont{\normalfont\rmfamily}}
\newtheorem{defi}[thm]{Definition}

\newtheorem{example}{Example}
\newtheorem{remark}{Remark}
\newenvironment{myexpcont}
{\addtocounter{example}{-1}\begin{example}{\textit{\textbf{{(continued)}}}}}
  {\end{example}}

\makeatletter
\def\BState{\State\hskip-\ALG@thistlm}
\title{An Algorithm for  Integer Least-squares with  Equality, Sparsity and Rank Constraints
}
\author{Arun Ayyar%
  \thanks{SantaFe Research Pvt Ltd., IIT Madras Research Park, Chennai--600113, India (\email{arun.ayyar@gmail.com}).}%
  \and
  Nirav P. Bhatt%
  \thanks{Robert Bosch Centre for Data Science and Artificial Intelligence, Indian Institute of Technology Madras, Chennai --600036, India
    (\email{niravbhatt@iitm.ac.in}).}
   }

\begin{document}
\maketitle

\begin{abstract}
In this work, we deal with rank-constrained integer least-squares optimization problems arising in low-rank matrix factorization related applications. We propose a solution for constrained integer least-squares problem subject to equality, sparsity, and rank constraints. The algorithm combines the  Fincke-Pohst enumeration  (or sphere decoding algorithm) with rank constraints and sparse solutions of Diophantine equations to arrive at an  optimal solution. The proposed approach consists of two steps as follows: (i) find the solution set for Diophantine equations arising from the linear and sparsity constraints, (ii) find the matrix which minimizes the integer least-squares objective and satisfying the rank constraints using the solution  set obtained in the step 1.  The proposed algorithm is illustrated using a simple example. Then,  we perform experiments to study the computational aspects of different steps of the proposed algorithm. 

\end{abstract}



\section{Introduction}
We consider constrained integer least-squares problems of the following form:
\begin{equation}
\begin{aligned}
& \underset{\Xbl}{\text{minimize}}
& &  \|\Ybl - \Gbl\Xbl \|^2   \\
& \text{subject to}
&&  \Xbl(i,j) \in \Scal \subset \Zcal,\, \forall i=1,\ldots,N,\,j=1,\ldots,L \\
&&& \boldsymbol{A}\Xbl^{\tr} = \zeros \\
&&& \rank{\Xbl}=N\\
&&& \|\Xbl^{\tr}(:,i)\|_0\leq K,\, i=1,\ldots,N
\end{aligned} \label{Eq:INMFProblem}
\end{equation}
 where  $\Ybl \in \Rcal^{M \times L}$, $\Gbl \in \Rcal^{M \times N}$,   $\Xbl \in \Scal^{N\times L}$ with $\Scal \subset \Zcal$, and $\boldsymbol{A} \in \Zcal^{P \times L} $. The equality constraints $\Abl\Xbl^{\tr} = \zeros$ are imposed on the rows of $\Xbl$. $\|\cdot\|^2$ denotes the square of  Frobenius norm. $||\cdot||_0$ is the zeroth norm of a given vector. The zeroth norm of a vector indicates the number of non-zero elements in a vector. The main motivation for investigating Problem~\ref{Eq:INMFProblem} is to solve subproblem arising in low-rank matrix factorization problems with integer constraints \cite{BerryBLPP07,KimP11,LeeH01,lin2011discrete}. For example, let us consider a matrix factorization problem in which a data matrix $\Ybl$ is decomposed into an integer matrix ($\Xbl$) and a non-negative matrix ($\Gbl)$ with the linear equality, rank and sparsity constraints  as follows \cite{Kim2008,Vavasis2009}:
\bea{ \Ybl \approx \Gbl \Xbl +\Ebl\label{Eq:integerFact1}}
where $\Ebl \in \Rcal^{M \times L}$ is an error matrix. 
The factorization can be achieved by minimizing a square of Frobenius norm of the error matrix as follows:
\begin{equation}
\begin{aligned}
& \underset{\Xbl,\, \Gbl}{\text{minimize}}
& &  \|\Ybl - \Gbl\Xbl \|^2   \\
& \text{subject to}
& &\Gbl \succeq \zeros\\
&& & \Xbl(i,j) \in \Scal \subset \Zcal,\, \forall i=1,\ldots,N,\,j=1,\ldots,L \\
&&& \Abf\Xbl^{\tr} = \zeros \\
&&& \rank{\Xbl}=N\\
&&& \|\Xbl^{\tr}(:,i)\|_0\leq K,\, i=1,\ldots,N
\end{aligned} \label{Eq:INMFProblem2}
\end{equation}
where $\succeq 0$ indicates each element of the matrix $\Gbl$ is greater or equal to zero. 
The objective function in \eqrefn{Eq:INMFProblem2} is a non-convex function, and it can be solved using alternating least-squares (ALS) framework.  Two  subproblems can be written  in the ALS framework as follows.
\begin{itemize}
\item For given $\Xbl$
\begin{equation}
\begin{aligned}
& \underset{\Gbl}{\text{minimize}}
& &  \|\Ybl - \Gbl\Xbl \| ^2  \\
& \text{subject to}
& &\Gbl \succeq \zeros\\
\end{aligned} \label{Eq:INMFProblemW}
\end{equation}
\item For given $\Gbl$
\begin{equation}
\begin{aligned}
& \underset{\Xbl}{\text{minimize}}
& &  \|\Ybl - \Gbl\Xbl \| ^2  \\
& \text{subject to}
& & \Xbl(i,j) \in \Zcal ,\, \forall i=1,\ldots,N,\,j=1,\ldots,L \\
&&& \Abl\Xbl^{\tr} = \zeros \\
&&& \rank{\Xbl}=N\\
&&& \|\Xbl^{\tr}(:,i)\|_0\leq K, i=1,\ldots,N
\end{aligned} \label{Eq:INMFProblemH}
\end{equation}
\end{itemize}
 
   The optimization problem  \ref{Eq:INMFProblemH} is a constrained integer least-squares (CILS) problem which is an NP-hard problem \cite{Schrijver00}. An efficient algorithm needs to be investigated for solving  this problem. In this work, we will propose an exact algorithm to solve Problem \ref{Eq:INMFProblemH}.
   
\subsection{Related Work in Literature} 
The unconstrained integer least-squares (ILS) problem is known as {\em closest point problem} in the lattice theory literature \cite{agrell2002closest}. The ILS problem also arises in several fields such as communications \cite{hassibi2005sphere}, global navigation systems\cite{teunissen1993least}, systems biology \cite{Bhatt2015} etc. Also, the ILS problem is a subset of nonlinear integer programming (NIP), or more specific, quadratic integer programming (QIP). Then,  methods to solve NIP or QIP can also be applied to solve the ILS problem.  Thus,  algorithms in the literature are based on approximation,  heuristics, and combinations of NIP and QIP for  solving the ILS problem \cite{borno2011reduction}. These algorithms consist of two steps: (i) Reduction, and (ii) search. In the first step, the columns of $\Gbl$ are orthogonalized as much as possible using lattice reduction methods such as Lenstra-Lenstra-Lov\'asz (LLL) or Korkine-Zolotareff (KZ)  reduction  \cite{borno2011reduction,lenstra1982factoring}. In the search step, the solution of the transformed ILS problem is obtained using Pohst enumeration (or also known as sphere decoder in communication literature) \cite{fincke1985improved, hassibi2002expected} or Schnorr-Euchner enumeration methods \cite{schnorr1994lattice}. Further, the constrained ILS problem with boxed or ellipsoid constraints have  been  solved in the literature of  communication field \cite{chang2009solving,chang2008solving}.

On the other hand, there have been efforts for developing efficient methods for  mixed integer nonlinear programming (MINLP), quadratic integer programming (QIP) or nonlinear integer programming (NIP)\cite{floudas1995nonlinear,hemmecke2009nonlinear}.  These methods can also be applied to  ILS problems. Several branch-and-bound  based algorithms have been developed to solve  MINLP, QIP and NIP problems in an efficient manner \cite{Billionnet2014,buchheim2012effective,buchheim2013exact,gupta1985branch}. A fast branch-and-bound algorithm for minimizing a convex QIP with convex constraints has been proposed \cite{buchheim2012effective}. The proposed algorithm by the authors allows to compute tighter lower bounds of the objective function by considering ellipsoidal under-estimators having the same continuous minimizer as the original objective function. The idea of ellipsoidal under-estimators has been extended to non-convex QIP with the box constraints \cite{buchheim2013exact}. This approach has been generalized in \cite{buchheim2015ellipsoid} using the under-estimators with strong rounding properties.   The proposed algorithm in \cite{buchheim2012effective,buchheim2013exact,buchheim2015ellipsoid} can also be applied to  ILS problems with boxed constraints in \cite{chang2009solving,chang2008solving}.

For non-convex QIP and MINLP problems, different approaches based on semi-definite relaxations have been proposed to obtain tighter lower bounds of the objective function \cite{buchheim2013semidefinite,Park2017}. The authors in \cite{buchheim2013semidefinite} proposed an algorithm in branch-and-bound approach by incorporating semi-definite relaxation for unconstrained non-convex MINLP problems. On the other hand, an approach to obtain lower and upper bounds on the objective function convex QIP using semi-definite relaxation is proposed in \cite{Park2017}. The algorithm proposed in \cite{Park2017} leads to a sub-optimal solution of the underlying problem. However, the computational experiments has shown that it provides near-optimal solutions for the larger size problem ($n=1000$). Moreover, Li et al. \cite{li2006convergent} and Saxena et al. \cite{saxena2010convex} have proposed algorithms based on contours and planes cuts for solving MINLP problems with inequality constraints. Recently, a fast branch-and-bound algorithms based on ellipsoidal relaxation for solving non-convex QIP with linear equality constraints has been proposed \cite{buchheim2015fast}. 

Most of the approaches in the literature can handle convex constrains or non-convex constraints like box, ellipsoid, equality, or inequality constraints. However, best of our knowledge, there is no algorithm which can solve Problem~\ref{Eq:INMFProblemH}.

\subsection{Contribution}
The main contribution of this paper is development of an algorithm for matrix integer-least squares problem in  \eqref{Eq:INMFProblemH}  with equality, sparsity and rank constraints by combining techniques from solution of a system of linear Diophantine equations and  modified Pohst enumeration (or modified sphere decoder). The proposed algorithm  consists of two steps:
\begin{enumerate}
	\item Find sparse solutions of a set of  Diophantine equations
	\item Find a solution matrix satisfying the rank constraints by combining solutions from the modified Pohst enumeration(or sphere  decoding) algorithm with the solution set obtained in the step 1.
\end{enumerate}
\section{Proposed Algorithms to Solve Problem~\ref{Eq:INMFProblemH}}\label{Algorithm} This section proposes a novel algorithm to solve the CILS problem \ref{Eq:INMFProblemH}. 
 There are three types of constraints in Problem~\ref{Eq:INMFProblemH}: (i) linear equality constraints, (ii) sparsity constraints, and (iii) rank constraint. The entries of $\boldsymbol{X}$ should be such that $\boldsymbol{X}^T$ is in the null space of $\boldsymbol{A}$ and each row of $\boldsymbol{X}$ can have atmost $K$ non--zero elements and the rank of $\boldsymbol{X}$ must be $N$. Further, the rank constraint can be checked for only after all the entries of $\boldsymbol{X}$ have been found. 
In order to find a solution to Problem~\ref{Eq:INMFProblemH}, we split the problem into two sub--problems based on the constraints.
\begin{enumerate}
\item The linear and sparsity constraints in Problem \ref{Eq:INMFProblemH} dictate that the rows of $\boldsymbol{X}$ must lie in the null space of $\boldsymbol{A}$ and the number of non-zero elements in the rows of $\boldsymbol{X}$ cannot be more than $K$. Then, the first problem is to find all the vectors that form this search space. Formally, Find the solution set $\mathcal{F}$  for
\begin{equation}
\mathcal{F}=\{\boldsymbol{x} \mid \boldsymbol{x}\in\mathcal{S}^L, \mathcal{S} \subset \mathcal{Z} \hspace{2mm} \& \hspace{2mm} \boldsymbol{A}\boldsymbol{x}=\boldsymbol{0}_{P \times 1} \hspace{2mm} \&  \hspace{2mm} ||\boldsymbol{x}||_0 \leq K \} \label{eq:prob1}
\end{equation} 

A set of  $N$ vectors from  $\mathcal{F}$  forms the rows of $\boldsymbol{X}$.

\item  Consider  $\boldsymbol{X}=[\boldsymbol{x}_1, \hspace{2mm} \boldsymbol{x}_2, \cdots, \boldsymbol{x}_L]_{N \times L}$ and $\boldsymbol{Y}=[\boldsymbol{y}_1, \hspace{2mm} \boldsymbol{y}_2, \cdots, \boldsymbol{y}_L]_{M \times L}$. Then, the objective function in \eqref{Eq:INMFProblemH} can be written as follows:
\begin{equation} 
\min_{\boldsymbol{X}\in S ^{N \times L}}||\boldsymbol{Y} - \boldsymbol{GX}||^2 = \sum_{j=1}^L \min_{\boldsymbol{x}_j\in \mathcal{S}_j}||\boldsymbol{y}_j - \boldsymbol{Gx}_j||^2\label{eq:prob2}
\end{equation}
The solution to each term in the summation in  \eqref{eq:prob2} is obtained by modified sphere decoding algorithm. Then, the next step is to pick $N$ vectors from $\mathcal{F}$ with help of the solutions for the columns of $\Xbl$ using the modified sphere decoding algorithm such that:
\begin{equation}
\text{Find }\boldsymbol{X} \text{ such that} \sum_{j=1}^L \min_{\boldsymbol{x}_j\in \mathcal{S}_j}||\boldsymbol{y}_j - \boldsymbol{Gx}_j||^2 \text{ and rank}(\boldsymbol{X})=N  \label{Eq:subproblem2}
\end{equation}

\end{enumerate}
In summary, \eqrefn{eq:prob1} provides a set candidates for the rows of $\Xbl$ and \eqrefn{eq:prob2} provides the solution set for the columns of $\Xbl$, while \eqrefn{Eq:subproblem2} finds an optimal solution of Problem~\ref{Eq:INMFProblemH} by combining the solutions of Eqs.~\eqref{eq:prob1}--\eqref{eq:prob2}.  The following example is used to illustrate each stage of the proposed algorithm in this section.\\
\begin{example}\label{exm1}
Consider the following matrices, and parameters
\begin{eqnarray}
\boldsymbol{Y} &=& \bMatsmall{0.5&3.7&-0.8&3.3&0.3&-3.5&-3.5\\1.8&5.8&-0.5&-0.4&-1.3&-2.7&-2.7\\-2.2&-3.1&2.6&0.5&-0.4&1.3&1.3\\0.8&3.5&-1.1&2.5&0.3&-3&-3}; \hspace{2mm} \boldsymbol{G}=\bMatsmall{0.5&0.3&3.5\\1.8&-1.3&2.7\\-2.2&-0.4&-1.3\\0.8&0.3&3} \nonumber\\
\boldsymbol{X}_{a}&=&\bMatsmall{1&1&-1&-1&0&0&0\\0&-1&-1&1&1&0&0\\0&1&0&1&0&-1&-1} \hspace{2mm} 
\boldsymbol{A}=\bMatsmall{8&2&10&0&12&2&0\\4&6&9&1&14&5&2\\2&0&1&1&0&1&0\\2&1&3&0&4&0&1}\label{eq:example1}
\end{eqnarray}
$\Scal=\{-1,0,1\}$, $K=4$, and $N=3$.
\end{example}

\subsection{Solution for a System of Linear Diophantine Equations}

In this section, we discuss the method to solve the problem described in \eqrefn{eq:prob1}. To solve $\Abl\Xbl^{\tr}=\zeros$ with sparsity constraint, we first transform  $\Abl$ in to a upper triangular matrix $\Hbl$ using the Hermite Normal form as follows \cite{Schrijver00}. 
\begin{defi}[Hermite Normal Form]
 For every $P \times L$ matrix $\boldsymbol{A}$ with integer entries, there exists a $P \times L$ matrix $\boldsymbol{H}$ with integer entries such that 
\begin{equation}
\boldsymbol{H}=\boldsymbol{UA} \hspace{3mm}\text{with}\hspace{3mm} \boldsymbol{U} \in GL_n(\mathcal{Z}) \label{Eq:hermiteNormalForm}
\end{equation}
where $\boldsymbol{U}$ is a $P \times P$-dimensional unimodular matrix and $\boldsymbol{H}$ is an upper triangular matrix. 
$\boldsymbol{H}$ is the matrix in Hermite normal form which can be obtained from $\boldsymbol{A}$  by elementary row operations. 
\end{defi}

Since $\boldsymbol{U}$ is a non-singular matrix in \eqrefn{Eq:hermiteNormalForm}, $\Abl\xbl=\zeros$ in \eqrefn{eq:prob1} can be written as
\begin{align}
\boldsymbol{UAx}&=\boldsymbol{0}  \\
\boldsymbol{Hx}&=\boldsymbol{0} \label{eqn:2}
\end{align}  
Since $\boldsymbol{H}$  is an upper triangular matrix, \eqrefn{eqn:2} can be written as
\begin{equation}
\begin{bmatrix} h_{1,1} & h_{1,2} & h_{1,3}& h_{1,4}&\cdots& h_{1,k}&\cdots & h_{1,L} \\ 0 & h_{2,2}& h_{2,3} & h_{2,4}&\cdots & h_{2,k} &\cdots & h_{2,L}\\ 0 & 0& 0 & h_{3,4}&\cdots & h_{3,k} &\cdots & h_{3,L}\\ \vdots &&&&&& \\ 0 & 0 & 0 & 0& \cdots & h_{P,k} & \cdots & h_{P,L} \end{bmatrix}\begin{bmatrix}x_1\\x_2\\x_3\\ \vdots \\ x_L
\end{bmatrix}=\begin{bmatrix}0\\0\\0\\ \vdots \\ 0\end{bmatrix} \label{eq:3}
\end{equation}
\eqrefn{eq:3} leads to the following set of equations
\begin{align}
h_{1,1}x_1 + h_{1,2}x_2+h_{1,3}x_3+\cdots+h_{1,L}x_L&=  0 \nonumber\\
h_{2,2}x_2+h_{2,3}x_3+\cdots+h_{2,L}x_L&=0\nonumber\\
&\vdots \nonumber\\
h_{P,k}x_k+\cdots+h_{P,L}x_L&=0 \label{eq:set_Dio}
\end{align}

The upper triangular structure of $\Hbl$ can be exploited in finding solutions for \eqrefn{eq:3}. Before explaining the methodology of finding the solution for a system of linear equations with constraints, we illustrate a method to find all solutions of a single linear Diophantine equation with $L_0$ norm constraint. The set of equations in \eqref{eq:set_Dio} can be written as follows:
 \bea{\hbl_i^{\tr}\xbl =0,\, i=1,\ldots,P}
 where $\xbl\in \mathcal{S}^{N_g}$ with $N_g\leq M$. Then, the objective is to 
 \begin{equation}
\text{ Find all }\boldsymbol{x} \in \mathcal{S}^{N_g} \text{ such that } \hbl_i^{\tr}\xbl = 0 \text{ and } ||\boldsymbol{x}||_0\leq K   \label{Eq:problem_16}\end{equation}
 Now, let us divide the variables $\xbl$ in to: (i) pivot variable ($x_{p}$) and (ii) free variables $(\xbl_{f})$. The free variables can take all the possible values from $\mathcal{S}$. The left most variable can be taken as the pivot variable and the rest variables can be taken as the free variables. Then, the values of $x_p$ in the $i$th equation can be computed as follows
 \bea{x_{p,i}=-\frac{\sum\limits_{l=1}^{N_g - 1}h_{i,l}x_{f,l}}{h_{p,i}} \label{Eq:pivotvariable}}
 
 Then, we need to find all possible solutions for $\xbl_f$.  Each element of $\xbl_f$ takes all the values in $\mathcal{S}$. For example, there are three elements in $\xbl_f$, i.e., $\xbl_f=[x_f(1),\, x_f(2),$
 $\, x_f(3)]^{T}$.  Then, the possible values for each variable are as follows: $\mathcal{S}_{x_{f}(1)}=\mathcal{S}_{x_{f}(2)}=\mathcal{S}_{x_{f}(3)}=\{-1,\,0,\,1\}$ $=\mathcal{S}$. Let us define the solution set till the $l$th variable as $\mathcal{S}_{\xbl_{f},l}$. Then, the set containing all the solutions for the  first variable, $\mathcal{S}_{x_f,1}=\mathcal{S}$.  Then, the set containing all the solutions till the  two variables variable $\xbl_{f,2}=[x_{f}(1),\, x_{f}(2)]^{T}$ can be obtained by a Cartesian product (denoted as $\times$) between $\mathcal{S}_{x_{f,1}}$ and the all solution $\mathcal{S}_{x_{f}(2)}$ as follows
 \begin{align}
     \mathcal{S}_{\xbl_f,2} = \mathcal{S}_{x_{f,1}} \times \mathcal{S}_{x_{f}(2)} = \mathcal{S} \times \mathcal{S}= \{(x_1,x_2) |  x_1 \in \mathcal{S}_{x_{f}(1)} \,\text{and}  \,x_2 \in \mathcal{S}_{x_{f}(2)} \}
     \label{Eq:two_state}
     \end{align}
Similarly, the solution set containing all three variables can be obtained by a Cartesian product between the previous stage solution set, $\mathcal{S}_{\xbl_f,2}$ and the all the possible solution for the third variable $x_{f}(3)$, $\mathcal{S}$ as follows:
 \begin{align}
     \mathcal{S}_{\xbl_{f,3}} = \mathcal{S}_{x_{f,2}} \times \mathcal{S}_{x_{f}(3)} = \mathcal{S}_{x_{f,2}} \times \mathcal{S}= \{(x_i) |  x_i \in \mathcal{S}_{x_{f}(i)},\, \forall i=1,2,3  \}
     \label{Eq:three_state}
     \end{align}
  \eqrefn{Eq:three_state} can be interpreted as follows.  The solution set for the three variables $\mathcal{S}_{\xbl_{f,3}}$ is a Cartesian product of the solution  set for the two variables $\mathcal{S}_{\xbl_{f,2}}$ and the set consisting all the possible values for $x_f(3)$, i.e.  $\mathcal{S}$. 
 Then, the solution set for the first  $l$ variables, $\xbl_{f,l}$, can be expressed as a Cartesian product of the solution set for the $l-1$ variables, $\mathcal{S}_{f,l-1}$ and the set consisting all possible values for the $l$ element of $\xbl_f$, $x_{f}(l)$,  $\mathcal{S}$ as follows:
 \begin{equation}
      \mathcal{S}_{\xbl_{f,l}} = \mathcal{S}_{x_{f,l-1}} \times \mathcal{S} \label{Eq:level}
 \end{equation}
 Since we are interested in the solutions that satisfy the sparsity constraint ($L_0-$norm) in   Problem \eqref{Eq:problem_16}, this constraint needs to be imposed on the set $\mathcal{S}_{\xbl_{f,l}}$ in \eqrefn{Eq:level} at each stage. To do so, let us define a function $f_s:D \rightarrow B$
 \begin{equation}
     f_s=\{\forall d\in D\,\,|\,\, \# (i|d_i\neq 0)\,\, \leq\,\, K\} \label{Eq:function_sparsity}
 \end{equation}
 The function in \eqref{Eq:function_sparsity}  imposes $L_0$-norm constraints on the set $D$ and find solutions satisfying the constraints. 
Then, the solution set till the $l$ variables  satisfying sparsity constraint ($L_0$-norm), $\hat{\mathcal{S}}_{\xbl_{f,l}}$, can be obtained as follows:
\begin{equation}
    \mathcal{\hat S}_{\xbl_{f,l}}=f_s(\mathcal{S}_{\xbl_{f,l}})=f_s(\mathcal{S}_{x_{f,l-1}} \times \mathcal{S})
    \label{Eq:sparse_solution_set}
\end{equation}
 Then, the complete solution set for  Problem \eqref{Eq:problem_16} can be obtained by  applying 
 \eqrefn{Eq:sparse_solution_set} till $l=N_g-1$.

 This solution set at any stage for \eqrefn{Eq:problem_16} can be stored using the concept of a rooted tree. Then, this rooted tree can be expanded to add more  variables by solving a set of equations in \eqrefn{eq:3}.
 Next, a methodology to obtain  the solution set for a system of linear Diophantine equations with $L_0$-norm constraint is described briefly. 
\begin{itemize}
\item Reduce the matrix $\boldsymbol{A}$ to its Hermite Normal form.
\item Start from the last equation.
\item Initialize the solution tree to the root node. 
\item For the free variables in the present equation, expand the existing solution tree. Check for $L_0$ norm constraint at each node.
\item For the pivot variable in the present equation, compute the pivot variable values using \eqrefn{Eq:pivotvariable} and also if the $L_0$ norm constraint is satisfied.
\item Retain only the final leaf nodes that satisfy both -  the current equation and the $L_0$ norm constraint.
\item Start from the solution tree for the previous equation and repeat the procedure for the next equation.
\item  The solution tree that remains after applying the method to all equations in \eqrefn{eq:3} is the final solution set $\mathcal{F}$ ans stored as tree T
\end{itemize}

The detailed method is given in Algorithm~\ref{algo:tree1}. Next, the elements of the algorithms are demonstrated using Example~\ref{exm1}.  

\begin{algorithm}
\caption{Algorithm to find solution to set of Diophantine equations with sparsity constraints  (Problem \ref{eq:prob1})\label{algo:tree1}}
\begin{algorithmic}[1]
\Require{ $K$, $\mathcal{S}$, $\boldsymbol{A} \in \Zcal^{P \times L} $}
\Ensure{The set $\mathcal{F}$ in terms of tree $T_{sol}$ containing all the solution vectors $\boldsymbol{x}$ such that $\boldsymbol{x}\in\mathcal{S}^{L}$ and $||\boldsymbol{x}||_0 \leq K$}
\Procedure{SolveSysOfLinDioEqn}{$\boldsymbol{A}$, $\mathcal{S}$, K}
\State $\boldsymbol{H}=\boldsymbol{UA}$ \Comment{Harmit Normal Form of $\boldsymbol{A}$}
\State $T=\{\}$
\State $N_{prev}=0$
\State $N$ = number of rows in $\boldsymbol{H}$
\For{$i=N$ to 1}
\State $I_d = $ Index of First non-zero element in $i^{th}$ row of $\boldsymbol{H}$ starting from left (Index of left most element is $0$)
\If{ $|T|$ ==0}
\State Add all elements of $\mathcal{S}$ as individual nodes to $T$
\Else
\State $N_{prev}$ = length of each node in $T$
\For{$k=1$ to $I_d-N_{prev}-1$}
\State $\boldsymbol{Z}$ = All the leaves in the Tree $T$
\State $N_z$ = Number of elements in $\boldsymbol{Z}$ 
\State  $T=\{\}$ \Comment{Empty the tree}
\For{$j=1$ to $N_z$}
\State$\boldsymbol{z} = \boldsymbol{Z}(j)$ \Comment{Assign $j^{th}$ vector in $\boldsymbol{Z}$ to $\boldsymbol{z}$}
\For {$l=1$ to $|\mathcal{S}|$}
\State$\boldsymbol{z} = [\boldsymbol{z}\hspace{2mm} \mathcal{S}(l)]$ \Comment{Append $l^{th}$ element in $\mathcal{S}$ to $\boldsymbol{z}$}
\If {$||\boldsymbol{z}||_0\leq K$}
 \State Add $\boldsymbol{z}$ as node to $T$  \Comment{Add the vector if it satisfies $L_0$ norm constraint}
\EndIf
\EndFor
\EndFor
\EndFor
\State $\boldsymbol{Z}$ = All the leaves in the Tree $T$
\State $N_z$ = Number of elements in $\boldsymbol{Z}$ 
\State  $T=\{\}$ \Comment{Empty the tree}
\State $\boldsymbol{h}_i=\boldsymbol{H}(i)$ \Comment{$i^{th}$ row of $\boldsymbol{H}$}
\For{$j=1$ to $N_z$}
\State$z_j = \boldsymbol{Z}(j)\cdot\boldsymbol{h}_i(I_d+1:end)/\boldsymbol{h}_i(I_d)$
\State $\boldsymbol{z} = [\boldsymbol{z}\hspace{2mm} z_j]$
\If {$z_j\in \mathcal{S}$ and $||\boldsymbol{z}||_0\leq K$}
 \State Add $\boldsymbol{z}$ as node to $T$  
\EndIf
\EndFor
\EndIf
\EndFor
\State \Return $T$
\EndProcedure
\end{algorithmic}
\end{algorithm}

\begin{myexpcont}
$\boldsymbol{A}$ can be factorized into its Hermite normal form as follows:
\begin{equation}
\underbrace{\left[\begin{smallmatrix}2& 0& 0& 2&-2&-8& 10\\0& 1& 0& 1& 0&-19& 21\\0 & 0& 1&-1& 2& 9&-10\\0& 0& 0& 0& 0&-18& 18\end{smallmatrix}\right]}_{\boldsymbol{H}} = \underbrace{\left[\begin{smallmatrix}-3&-1& 3& 12\\-6&-2& 3& 25\\3& 1&-2&-12\\-5&-2& 2& 22\end{smallmatrix}\right]}_{\boldsymbol{U}}\underbrace{\left[\begin{smallmatrix}
8&2&10&0&12&2&0\\4&6&9&1&14&5&2\\2&0&1&1&0&1&0\\2&1&3&0&4&0&1\end{smallmatrix}\right]}_{\boldsymbol{A}} \label{eq:hnf1}
\end{equation}
The system of equations $\boldsymbol{Ax}=\boldsymbol{0}$ can thus be written as
\begin{align}
2x_1+ 2x_4-2x_5-8x_6+ 10x_7&=0 \nonumber\\
x_2+x_4-19x_6+21x_7&=0 \nonumber\\
x_3-x_4+ 2x_5+ 9x_6-10x_7&=0 \label{eq:setofeqns} \\
-18x_6+18x_7&=0 \nonumber
\end{align}
The last equation in \eqref{eq:setofeqns} is considered to demonstrate building of the solution tree. The problem to be solved is as follows:
\begin{equation}
-18x_6+18x_7=0 \label{eq:lasteqn}
\end{equation}
such that $||\boldsymbol{x}||_0\leq4$ and $\boldsymbol{x} \in \{-1,0,1\}^2$, where $\boldsymbol{x}=\begin{bmatrix}x_6 & x_7 \end{bmatrix}$. We need to find  all possible solutions of \eqrefn{eq:lasteqn}.
In \eqref{eq:lasteqn}, $x_6$ is the pivot variable and $x_7$ is the free variable. Hence, $x_7$ can take all the values from $\mathcal{S}$ while the values of  $x_6$ depend on $x_7$.  We start with an empty tree. At each stage, we consider the values that a particular variable can take. For example,  in first stage, $x_7$ can take all possible values from $\mathcal{S}$. Hence, to the root node, we add 3 nodes one each for the values $-1$, $0$, and $1$. At each subsequent stage, we concatenate the value taken by a variable to the value present at its parent node. In the next stage, we consider the variable $x_6$.  $x_6$  can be computed using \eqrefn{Eq:pivotvariable}. Then, the solution set for \eqrefn{eq:lasteqn} is shown in Figure.~\ref{fig:tree15} in terms of a rooted tree. 


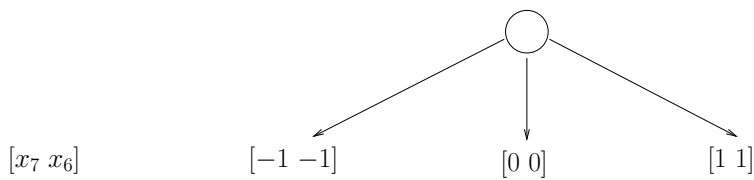
\begin{figure}[h]
\centering
\scalebox{0.5}
{\input{reverse_sol_tree2.pstex_t}}\vspace{-1.5cm}
 \caption{Solution tree to \eqrefn{eq:lasteqn}}
\label{fig:tree15}
\end{figure}


Then,  by applying Algorithm~\ref{algo:tree1}, the solution tree for  \eqrefn{eq:setofeqns} is given in Fig.~\ref{fig:extree1}. The intermediate solutions are also depicted.  

\begin{figure}[h]
\centering
\scalebox{0.45}
{\input{reverse_sol_tree4.pstex_t}}\vspace{-0.25cm}
 \caption{Solution set $\mathcal{F}$ for Example~\ref{exm1}.}
\label{fig:extree1}
\end{figure}
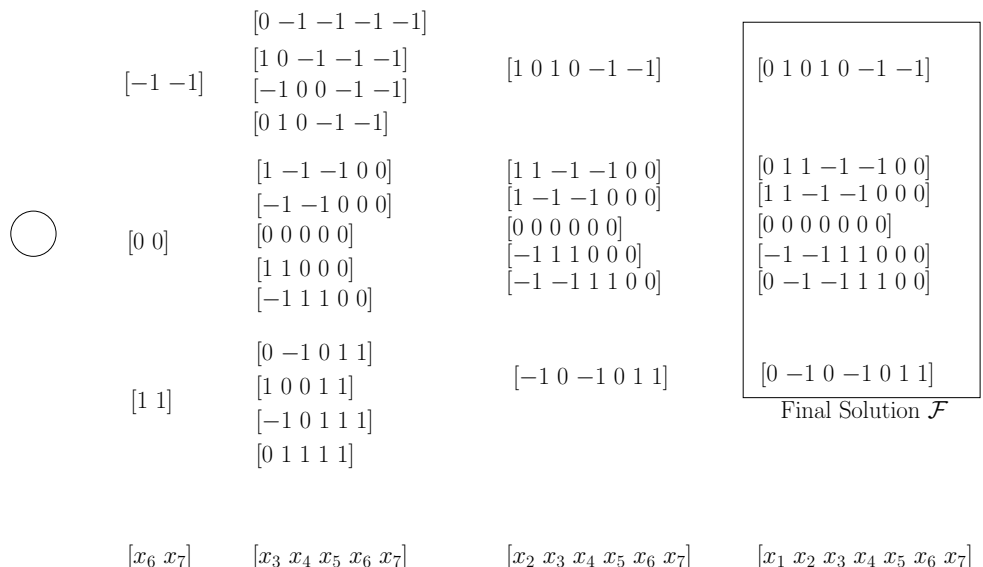
There are seven solutions that satisfy $\Abl\Xbl^{T}=\zeros$ with $K=4$ which are shown in the box and labelled as {\em the final solution} in Figure~\ref{fig:extree1}. The solution set $\mathcal{F}$ contains all possible rows of $\Xbl$.
\end{myexpcont}


\subsection{Solution to Constrained Integer least-squares}
In the previous section, we outlined the method to find the solution to the sub-problem~\eqref{eq:prob1}. This section proposes an approach   to solve the posed CILS problem \ref{Eq:INMFProblemH} using the solution set $\mathcal{F}$. To solve the proposed problem, we first revisit the algorithm to solve of  Problem~\ref{eq:prob2} using the modified sphere decoding method.  The  individual terms within the summation in \eqrefn{eq:prob2} can be formulated as the following problem. Given $\mathcal{S}\subset \mathcal{Z}^N$,  $\boldsymbol{G}\in\mathcal{R}^{M \times N}$ and $\boldsymbol{y} \in R^{M}$, 
 \begin{equation}
  \min_{\boldsymbol{x_j}\in \mathcal{S}}||\boldsymbol{y}_j - \boldsymbol{Gx_j}||^2. \label{eq:ILS}
 \end{equation} 
The problem presented in \eqref{eq:ILS} is the classical {\em integer least--squares} problem. Since the entries of $\boldsymbol{x}_j$ are restricted to $\mathcal{S} \in \mathcal{Z}^N$, the search space is a finite subset of the infinite ``rectangular" lattice $\mathcal{Z}^N$. When multiplied by a matrix $\boldsymbol{G}$, $\boldsymbol{Gx_j}$ spans a ``skewed lattice". Hence, given a ``skewed lattice" $\boldsymbol{Gx_j}$ and a vector $\boldsymbol{y}_j$, the integer least--squares problem is to find the lattice point closest to $\boldsymbol{y}_j$ in the Euclidean sense. However, in contrast to classical ILS problem \eqref{eq:ILS}, it is possible that each element of $\xbl_j$ can belong to a different set in this work, i.e., $\xbl_j(i)\in\mathcal{S}_{i,j} \subseteq \mathcal{S},\,\forall i=1,\ldots,N$.  By defining a set $\mathcal{S}_j=\{\mathcal{S}_{1,j},\,\mathcal{S}_{2,j},\,\ldots \mathcal{S}_{N,j}\}$, Problem~\ref{eq:ILS} can be written as follows
 \begin{equation}
  \min_{\boldsymbol{x}_j\in \mathcal{S}_j}||\boldsymbol{y}_j - \boldsymbol{Gx}_j||^2. \label{eq:ILS2}
 \end{equation} 
Then, Problem~\eqref{eq:ILS2} can be solved by the existing algorithms for solving the classical ILS problem in \eqrefn{eq:ILS}.
We can exploit the nature of the problem by modifying the existing algorithm to include the information of the set $\mathcal{S}_j$.  This helps us to reduce computational complexity. Hence, we will  modify the existing  sphere decoding algorithm in the current work \cite{fincke1985improved,hassibi2005sphere}.  
Then, the existing sphere decoding algorithm can be modified by  changing   the bounds on each element $i$ of $\xbl_j(i)$ as presented in Algorithm~\ref{Algo:spheredecoding}.

\begin{algorithm}
\caption{Modified sphere decoding algorithm to solve Problem \eqref{eq:ILS2}}\label{Algo:spheredecoding}
\begin{algorithmic}[1]
\Require{$\ybl$, $\Gbl$, Search radius $d$, Set $\mathcal{S}$}
\Ensure{Solution vectors $\boldsymbol{x}$}
\Procedure{SphereDecoder}{$\mathcal{S}$,$\ybl$,$\boldsymbol{G}$,$d$}
\State $\Gbl= [\Qbl_1,\,\Qbl_2]\bMat{\Rbl \\\zeros}$\Comment{ QR decomposition of $\Gbl$}
\State Compute Hermitian transpose  $\Qbl_1^*$ and $\Qbl_2^*$ of the sub-matrices $\Qbl_1$ and $\Qbl_2$, and $\zbl=\Qbl^*_1\ybl$.
\State Set $i=N$, $d'^2_N =d^2-\|\Qbl_2^*\ybl\|$, $z_{N|N+1}=z_N$.
\State Set $u_b(x_i)=\min(\max(\mathcal{S}_i),\lfloor\frac{d'_i+z_{i|i+1}}{r_{i,i}}\rfloor$),\,
\State $x_i = \max(\min(\mathcal{S}_i),\lceil \frac{-d'_i+z_{i|i+1}}{r_{i,i}}\rceil) -1 $
\Comment{Bounds for $x_i$}
\State $x_i=x_i+1$, 
\If  {$x_i \leq u_b(x_i)$}, go to 12 
\Else\,\,  go to 9
\EndIf
\State $i=i+1$ \Comment{Increase $i$}
\If {$i=N+1$}, terminate algorithm
\Else\,\,  go to 7
\EndIf
\If {$i=1$}, go to 14 \Comment{Decrease $i$}
\Else \,\, $i=i-1$, $z_{i|i+1} = z_i -\sum_{j=i+1}^{N}r_{i,j}x_j$, $d'^2_{i} =d'^2_{i+1} - (z_{i+1|i+2}-r_{i+1,i+1}x_{i+1})^2$ and go to 5
\EndIf
\State Solution found. Save $\xbl$, $d'^2_N -d'^2_1+ (z_1 -r_{1,1}x_1)^2$, and go to 7.
\State \Return $\xbl$
\EndProcedure
\end{algorithmic}
\end{algorithm}

They key concept that sphere decoding uses to solve the problem is a parameter $d$ which is the radius of the sphere within which the algorithm will search for $\boldsymbol{x}_j$. Assuming the centre of the sphere to be the given vector $\boldsymbol{y}_j$, the algorithm searches for the lattice points of $\boldsymbol{Gx_j}$ which lies within a sphere of radius $d$. It can be seen that, the closest lattice point inside the sphere will also be the closest lattice point for the whole lattice. If a point is found within the sphere, then the vector is returned, else the radius of the sphere is increased and the search is performed all over again. The details of  the sphere decoding algorithm  and its complexity can be found in \cite{fincke1983procedure}.

There is a rank constraint  in the sub-problem~\ref{Eq:subproblem2}. Since the closest lattice point inside the sphere for each column does not guarantee that $\Xbl$ is of rank $N$.  The output of algorithm~\ref{Algo:spheredecoding} provides all lattice points inside the radius of sphere $d$ which will be used to combine with the set $\mathcal{F}$

\begin{remark}
The  Sphere decoding returns candidate solutions for each of the columns of $\boldsymbol{X}$. 
As can be seen from \eqref{eq:prob2}, we need to run the  sphere decoder algorithm for all the $L$ columns of $\boldsymbol{X}$. The inputs to the sphere decoding for the solution of the $j${th} column of $\boldsymbol{X}$ will be as seen from  \eqrefn{eq:ILS2} are $\mathcal{S}_j^N$, $\boldsymbol{G}$ and $j^{th}$ column of $\boldsymbol{Y}$, i.e. $\boldsymbol{y}_j$. 
The rows of $\boldsymbol{X}$ are chosen from the solution set $\mathcal{F}$ of  \eqrefn{eq:prob1}. Hence, the search space for the sphere decoding need not be chosen as $\mathcal{S}^N$  for all the $L$ runs. The search space $\mathcal{S}_j$ for the $j^{th}$ run will be decided by the $i^{th}$ entry of the solutions to Problem \ref{eq:prob1}. This will help in reducing the overall search space while solving Problem~\ref{eq:ILS2} using the modified sphere decoding algorithm.\\
\end{remark}
As explained earlier, the candidates for the rows of $\boldsymbol{X}$, $\xbl_R$, are generated from the solutions to Problem~\ref{eq:prob1}, i.e. $\mathcal{F}$, and the column sets of $\boldsymbol{X}$, $\xbl_C$ are generated by Algorithm~\ref{Algo:spheredecoding}. Then, the objective is to  find the actual entries of $\boldsymbol{X}$ using the solutions to Eqs.~\eqref{eq:prob1} and \eqref{eq:prob2} such that the rank constraint is also satisfied. To achieve this objective, $\boldsymbol{X}$ has to be constructed from the solutions of rows  and columns candidates, and then check for the rank condition.  Let us define the solution set for $\xbl_C$ is $\mathcal{S}_C$ with $|\mathcal{S}_C|=L$ where $|\cdot|$ indicates the cardinality of a given set. Further, the solution set for $i$th row of $\Xbl$ at the computation of the $j$th column, $\mathcal{F}_{i,j},\,\,i=1,2,\ldots,N,\,\, j=1,2,\ldots,L$, is given by 
\begin{equation}
    \mathcal{F}_{i,j}=\{\xbl_R\,|\, \xbl_R(j)=\xbl_{C,j}(i),\,\forall \xbl_R \in \mathcal{F},\,\xbl_C \in \mathcal{S}_C\} \label{Eq:Set_Consistence}
\end{equation}
where $\xbl_R(j)$ is the $j$th element of the row vector $\xbl_R$ and $\xbl_{C,j}(i)$ indicates the $i$th element of the $j$th column in the ordered set $\mathcal{S}_C$. Then,  the sets $\mathcal{F}_{i,j}$ can be used to create $\Xbl$ matrices, and the rank of $\Xbl$ matrices can be checked.  One optimal way to implement this solution strategy is to incrementally solve problem for each column $\xbl_{C,j}$ and then refine the elements of  $\mathcal{F}_{i,j}$ at  $j$th iteration. Next, we describe briefly this approach. Note that $\mathcal{F}_{i,0}=\mathcal{F},\,\forall i=1,\ldots,N$.

Since the solution sets to \eqrefn{eq:prob1} are contained in the solution trees and there are $N$ rows in $\boldsymbol{X}$,  each of the $N$ rows is chosen from  the solution set $\mathcal{F}_{i,j}$. As a first step, we repeat the solution tree $N$ times, once for each row. Once Algorithm~\ref{Algo:spheredecoding} finds  the solution for a given column, the solution tree for each row will be pruned in order to find the candidate vector for each of the rows of $\boldsymbol{X}$ using the set definition in \eqrefn{Eq:Set_Consistence}. 
The method of finding the entries of $\boldsymbol{X}$ is briefly described as  follows:
\begin{itemize}
\item Find the solution tree containing solutions of the $\xbl_R$ to Problem~\ref{eq:prob1} using Algorithm~\ref{algo:tree1}. 
\item Repeat the solution tree $N$ times. 
\item Assume a value for $d$, the radius in Algorithm~\ref{Algo:spheredecoding}.
\item In the $i${th} solution tree (i.e. for $i${th} row of $\boldsymbol{X}$), denote the set of unique entries of the $j${th} column of the solutions present in the tree as $\mathcal{S}_{j}$. Define  $\mathcal{S}_j = [\mathcal{S}_{1,j} \hspace{2mm} \mathcal{S}_{2,j} \cdots \mathcal{S}_{N,j}]^T$. $\mathcal{S}_j$ is the input to Algorithm~\ref{Algo:spheredecoding} in the $j${th} run.
\item Algorithm~\ref{Algo:spheredecoding} for the $j^{th}$ column will output all candidates which are within the sphere radius. Denote the output of Algorithm~\ref{Algo:spheredecoding} as $\boldsymbol{Z}_j$. 
\item Run the Algorithm~\ref{Algo:spheredecoding} for all the columns.
\item For each of the columns, from $\boldsymbol{Z}_j$ choose the vector in the ascending order of distance which has not been considered previously. Denote this vector as  $\boldsymbol{{x}}_{C,j} = [{x}_{C,j}(1) \hspace{2mm} {x}_{C,j}(2) \ldots {x}_{C,j}(N)]^T$.
\item In the $i${th} solution tree, remove all the vectors whose $j${th} element is not same ${x}_{C,j}(i)$. This is equivalent to refining the solution sets $\mathcal{F}_{i,j}$. If only one vector is left in the tree do not perform the pruning step for the $i${th} tree.
\item Repeat this procedure for all the columns of $\boldsymbol{X}$.
\item The $i${th}  solution tree contains the candidate solutions for the $i${th} row of $\boldsymbol{X}$.
\item Form the matrix $\boldsymbol{X}$ and check for the rank condition. If the rank condition is satisfied, declare the derived matrix as the solution, else try out other vectors from $\boldsymbol{Z}_j$ which have not yet been considered.
\item If all the vectors in $\boldsymbol{Z}_j$ are exhausted and still the rank condition is not satisfied, then reduce $j$ by 1 and re-run the procedure.
\item If $j$ becomes zeros and still no solution is found, then increase $d$ and repeat the procedure till one $\boldsymbol{X}$ is found that satisfies the rank condition.
\end{itemize}
The complete algorithm for finding the solution to Problem~\ref{Eq:INMFProblemH} can be given in Algorithm~\ref{CompleteAlgo}. 

\begin{algorithm}[H]
\caption{Algorithm to solve Problem \eqref{Eq:INMFProblemH}}\label{CompleteAlgo}
\begin{algorithmic}[1]
\Require{$\boldsymbol{G}\in\mathcal{R}^{M \times N}$, $\boldsymbol{Y} \in \mathcal{R}^{M \times L}$ and $\boldsymbol{A} \in \mathcal{Z}^{P \times L}$, $\mathcal{S}$, $K$}
\Ensure{Matrix $\boldsymbol{X}\in S ^{N \times L}$ such that $\boldsymbol{AX}^T=\boldsymbol{0}$, $||\boldsymbol{X}^T(i)||_0 \leq K$, rank($\boldsymbol{X}$)$=N$ and $||\boldsymbol{Y}-\boldsymbol{GX}||^2$ is minimum} 
\Procedure{SolveSPDecWCons}{$\boldsymbol{G}$, $\boldsymbol{Y}$, $\boldsymbol{A}$ , $\mathcal{S}$, $K$}
\State $T$={SolveSysOfLinDioEqn}($\boldsymbol{A}$, $\mathcal{S}$, $K$)
\State $\boldsymbol{X}=\boldsymbol{0}$
\State Discard the vectors in $\boldsymbol{Z}_j$ that are already present in $\boldsymbol{Z}_{p,j}$

\State  Repeat the solution tree $T$, $N$ times, once for each row. 
\State Connect the root nodes of the  repeated trees to a common root node. 
\State Denote this tree as $T_N$
\While{rank($\boldsymbol{X}$)$< N$}
\State Set $\boldsymbol{Z}_{p,1},\cdots,\boldsymbol{Z}_{p,L}=\boldsymbol{0}$
\State $j=1$ , $J=L$,$M=0$
\While{$j \leq L$}
	\For {$i=1$ to $N$}
	\State $\mathcal{S}_{i,j}=$ Set of unique elements at $j^{th}$ position in  the $i^{th}$ branch  vectors of  $T_N$
	\EndFor
\State$\mathcal{S}_j=\{\mathcal{S}_{1,j},\mathcal{S}_{2,j},\cdots,\mathcal{S}_{N,j}\}$
\State$\boldsymbol{Z}_{j}=$SphereDecoder($\mathcal{S}_j$,$\boldsymbol{Y}(:,j)$,$\boldsymbol{G}$,$d_{sp}$)

\State $\boldsymbol{x}_{C,j}= $ Vector from $\boldsymbol{Z}_j$ which has not be considered previously in the ascending order of distance.
\For{$k=1$ to $N$}		\If{Number of vectors in $k^{th}$ branch of $T_N$ $ >1$}
		\State In $k^{th}$ branch of $T_N$, remove all vectors which dont have $\boldsymbol{x}_{C,j}(k)$ in the $j^{th}$ position
		\EndIf
	\EndFor
	
	\If{$j==L$}
	\State Create a matrix $\boldsymbol{X}$ containing one vector from each branch in $T_N$

	\If  {rank($\boldsymbol{X}$)$!=N$}
\State	$n$=Number of vectors in $\boldsymbol{Z}_j$
    \State$M=M+1$
    \If{$M==n$}
    \State $J=J-1$, $M=0$
    \EndIf
    
    \If{$J==0$}
    \State $j=1$,$J=L$
    \State $d_{sp}=d_{sp}+1$
	\For{$k=1$ to $L$}
\State	$\boldsymbol{Z}_{p,k}=\boldsymbol{Z}_k$
	\EndFor
	\Else
	\State $j=J$
	\State Go To 16
    \EndIf

	\EndIf
	\EndIf
	\State $j=j+1$
	\EndWhile
	
\EndWhile
\Return $\boldsymbol{X}$
\EndProcedure
\end{algorithmic}

\end{algorithm}

\begin{remark}
Note that if the equality constraints is of form $\Abl\Xbl=0$ instead of $\Abl\Xbl^{\tr}=0$, the construction of $\Xbl$ is straightforward using the set $\mathcal{F}$.
\end{remark}
Next, Algorithm~\ref{CompleteAlgo}
 is explained using  Example~\ref{exm1}.
\begin{myexpcont}

There are three rows in $\boldsymbol{X}$. Therefore, the solution set $\mathcal{F}$ in Fig.~\ref{fig:extree1} is replicated three times as shown in Fig.~\ref{fig:soltree1} one solution tree for each row.
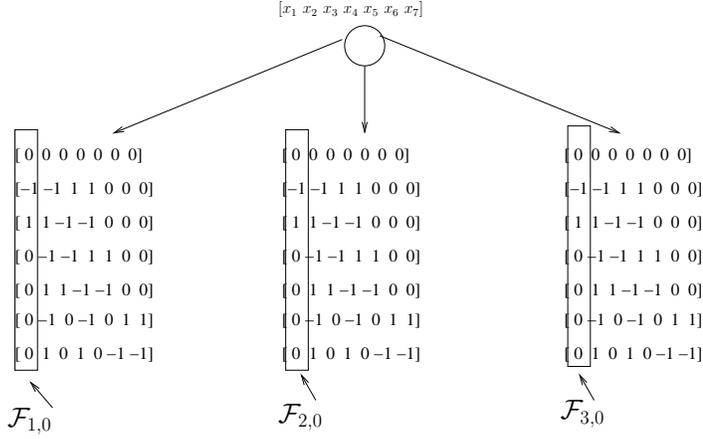
\begin{figure}[h]
\centering
\scalebox{0.5}
{\input{sol_tree1.pstex_t}}
 \caption{Solution tree  \eqref{eq:prob1} repeated $3$ times for each rows, and the corresponding the solution sets are denoted as $\mathcal{F}_{i,0},\,\forall i=1,2,3$}
\label{fig:soltree1}
\end{figure}
It can be seen that in the first tree set of unique entries for the  first column is $\mathcal{S}_{i,1}=\{-1,0,1\}$. The set of unique entries for all the $3$ rows is the same. Hence for the first column $\mathcal{S}_1= \{-1,0,1\}^3$.

Using the radius  $d=0.5$, the solution given Algorithm~\ref{Algo:spheredecoding} for the first column ($j=1$) of $\boldsymbol{X}$ is $\boldsymbol{{x}}_{C,1}=\begin{bmatrix}1 &0&0\end{bmatrix}^{T}$.
This means that the candidate vectors for the rows $R1$, $R2$, and $R3$ must have $1$, $0$ and $0$ in the first position. Then, the  new sets  $\mathcal{F}_{i,1},\, i=1,2,3$ are obtained by deleting the candidates in $\mathcal{F}_{i,0},\, i=1,2,3$ which do not contain the specified elements in $\xbl_{C,1}$ on the first position for the corresponding rows.  The resulting solution tree is shown in Fig. \ref{fig:soltree2}.
\begin{figure}[h]
\centering
\scalebox{0.5}
{\input{sol_tree2.pstex_t}}
 \caption{Solution tree for \eqref{eq:example1} and  pruning operation using  the solution for $\xbl_{C,1}$}
\label{fig:soltree2}
\end{figure}
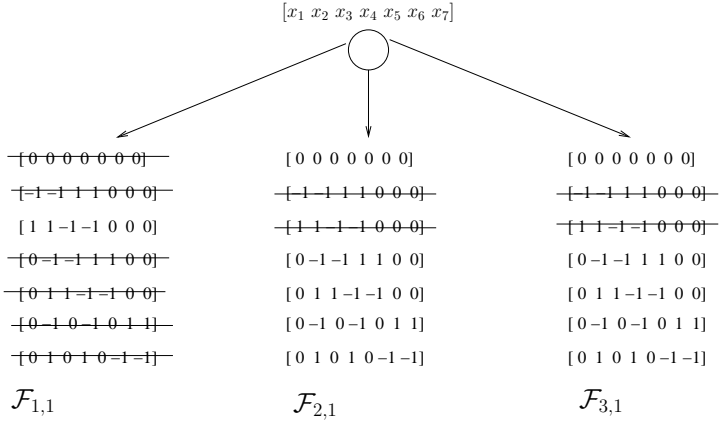
As seen from Fig. \ref{fig:soltree2}, there is only one candidate vector for the $\xbl_{R,1}$ of $\boldsymbol{X}$, while there are five candidate vectors for the rows   $\xbl_{R,2}$  and  $\xbl_{R,3}$. 
 As it can be seen from  Fig.~\ref{fig:soltree2}, the set for each element of  $\xbl_{C,2}$ can be re-defined as: $\mathcal{S}_2=\{\mathcal{S}_{1,2},\,\mathcal{S}_{2,2},\, \mathcal{S}_{3,2}\}$, $\mathcal{S}_{1,2}=\{1\}$, $\mathcal{S}_{2,2}=\{-1,0,1\}$ and  $\mathcal{S}_{3,2}=\{-1,0,1\}$. With these sets and  the solutions from Algorithm~\ref{Algo:spheredecoding} with $d=0.5$, we get the solution set for the second column ($j=2$),  $\boldsymbol{{x}}_{C,2}=\begin{bmatrix}1&-1&1\end{bmatrix}^T$. This means that  we need to eliminate all the vectors that do not  have $-1$ and $1$, respectively in their second position from the candidate solution sets $\mathcal{F}_{2,1}$ and $\mathcal{F}_{3,1}$.  The resulting tree is shown in Fig. \ref{fig:soltree3}.

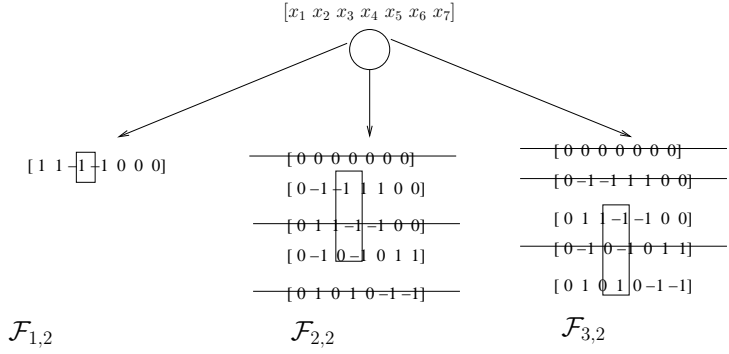
\begin{figure}[h]
\centering
\scalebox{0.5}
{\input{sol_tree3.pstex_t}}
 \caption{Solution tree for \eqref{eq:example1} and    pruning operation using the solution for $\xbl_{C,2}$}
\label{fig:soltree3}
\end{figure}
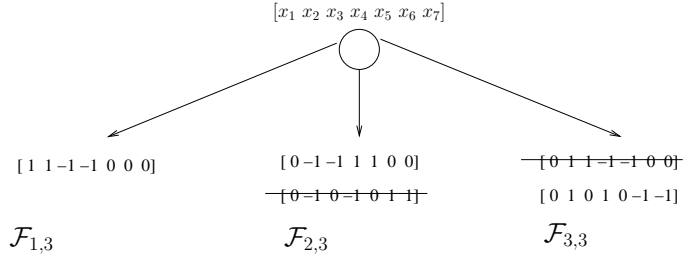
\begin{figure}[h]
\centering
\scalebox{0.5}
{\input{sol_tree4.pstex_t}}
 \caption{Solution tree for \eqref{eq:example1} and  pruning operation using the solution $\xbl_{C,3}$}
\label{fig:soltree4}
\end{figure}
Constructing the set of all allowable entries for $\xbl_{C,3}$ of $\boldsymbol{X}$, we have $\mathcal{S}_{1,3}=\{-1\}$, $\mathcal{S}_{2,3}=\{-1,0\}$ and $\mathcal{S}_{3,3}=\{0,1\}$.
Applying Algorithm~\ref{Algo:spheredecoding} to the third column ($j=3$) of $\boldsymbol{Y}$ and the corresponding set $\mathcal{S}_3$, we have $\boldsymbol{\hat{x}}_{C,3}=\begin{bmatrix}-1&-1&0  \end{bmatrix}^T$. Eliminating candidate vectors for the rows $R2$ and $R3$ that do not contain $\{-1\}$ and $\{0\},$ respectively, the solution tree is shown in Fig.~\ref{fig:soltree4}.  After pruning with $\xbl_{C,3}$, it  can be seen from Fig.~\ref{fig:soltree4} that the sets $\mathcal{F}_{i,3},\,\forall i=1,2,3$ contain only one element.

 Hence,  we cannot search any further since there are no other alternative solutions. Hence,  $\Xbl$ can be constructed as follows
\begin{equation}
\boldsymbol{{X}}=\begin{bmatrix}
1&1&-1&-1&0&0&0\\0&-1&-1&1&1&0&0\\0&1&0&1&0&-1&-1
\end{bmatrix}
\end{equation}
Since $\rank{\Xbl}=3$, the solution is optimal $\Xbl^*$ 
which is same as the actual solution $\boldsymbol{X}_a$ as given in \eqref{eq:example1}. Note that if $\rank{\Xbl}\neq 3$, then $d$  has to be increased in Algorithm~\ref{Algo:spheredecoding} and repeat the procedure. 

\end{myexpcont}

\section{Special case: Integer least--squares  with the linear equality and sparsity constraints}
In Section~\ref{Algorithm}, an algorithm to solve Problem~\ref{Eq:INMFProblemH} is proposed. The main advantage of the proposed algorithm is that the two parts solution to the the optimization problem \eqref{Eq:INMFProblemH}. Using this two-parts algorithm, we can also solve integer least--squares with the linear inequality and sparsity constraints. In this case, the optimization problem is:
\begin{equation}
\begin{aligned}
& \underset{\xbl}{\text{minimize}}
& &  \|\ybl - \Gbl\xbl \|^2   \\
& \text{subject to}
&&  x_i \in \Scal \subset \Zcal,\, \forall i=1,\ldots,N \\
&&& \boldsymbol{A}\xbl = \bbf\\
&&& \|\xbl\|_0\leq K.\, 
\end{aligned} \label{Eq:INMFProblemmodified1}
\end{equation}
The algorithm for the problem \eqref{Eq:INMFProblemmodified1} can be obtained by putting together the solutions of the sub-problems of the problem~\eqref{Eq:INMFProblem} as follows:
\begin{enumerate}
\item Find all the solutions of the  Diophantine Equations with the sparsity constraints. Formally, we would like to find a solution set  $\mathcal{F}$ in \eqref{eq:prob1}.
\item Then, solve the integer least--squares problem of the following form:
\begin{equation}
\text{Find }\boldsymbol{x} \text{ such that}  \min_{\boldsymbol{x}\in \mathcal{F}}||\boldsymbol{y} - \boldsymbol{Gx}||^2. \label{Eq:ILS_F}
\end{equation}
In contrast to the standard integer--least squares problem in \eqref{eq:ILS}, the solution space is constrained by the solution set $\mathcal{F}$ in  \eqrefn{Eq:ILS_F}. 
\end{enumerate}
An algorithm to solve the problem \eqref{Eq:INMFProblemmodified1} is described next using the solution strategies  in the last section.

\begin{algorithm}[H]
\caption{Algorithm to solve \eqref{Eq:INMFProblemmodified1}}
\begin{algorithmic}[1]
\Require{$\boldsymbol{G}\in\mathcal{R}^{M \times N}$, $\boldsymbol{y} \in \mathcal{R}^{M \times 1}$ and $\boldsymbol{A} \in \mathcal{Z}^{P \times L}$, $\mathcal{S}$, $K$}
\Ensure{Vector $\boldsymbol{x}\in S ^{N \times 1}$ such that $\boldsymbol{AX}=\boldsymbol{0}$, $||\boldsymbol{x}||_0 \leq K$, and $||\boldsymbol{y}-\boldsymbol{Gx}||^2$ is minimum} 
\Procedure{SolveILSEQ}{$\boldsymbol{G}$, $\boldsymbol{y}$, $\boldsymbol{A}$ , $\mathcal{S}$, $K$}
\State $\mathcal{F}$={SolveSysOfLinDioEqn}($\boldsymbol{A}$, $\mathcal{S}$, $K$)
\State$\hat{\boldsymbol{x}}=$SphereDecoder($\mathcal{S}_i$,$\boldsymbol{y}$,$\boldsymbol{G}$,$d_{sp}$)
\State $\dbl_x=\zeros$
\For{$i=1$ to $|\mathcal{F}|$}
\State Compute $\dbl_x(i) =\|\hat\xbl - \xbl_i \|^2$ where $\xbl_i$ is the $i$th element in the set $\mathcal{F}$
\EndFor
\State Find the minimum value in the $\dbl_x$ and the corresponding index $f$
\State $\xbl_{min}$ = $f$th element of the set $\mathcal{F}$

			 
	
\Return $\boldsymbol{x}_{min}$
\EndProcedure
\end{algorithmic}

\end{algorithm}
\section{Numerical Experiments}
This section  reports computational results obtained by applying the proposed algorithm. We consider simulation examples of different sizes and rank of $\Xbl$. For each simulation example, the data matrix $\Ybl$ were generated based on the model in  \eqrefn{Eq:integerFact1}. The matrix $\Gbl$  is generated as follows. The random samples of desired dimension ($\Gbl$) are generated   under the assumption of the Gaussian distribution with zero mean and the unity standard deviation. For obtaining positive matrix, the absolute of the generated matrix is taken.  Then, the matrix $\Xbl$  are generated by picking  such that the constraints are satisfied. We restrict each $X(i,j)$ element in $\mathcal{S}=\{-2,-1,0,1, 2\}$.   The sparsity $K=4$ is taken.   Each element of matrix $\Ebl$, $\Ebl(i,j)$ is considered to be normally distributed with zero mean and standard deviation $\sigma=0.2$.  The data matrix, $\Ybl$ is obtained by the multiplying the generated  $\Gbl$ and $\Xbl$ and adding  $\Ebl$ of the appropriate size. The equality constraints are generated using the matrix $\Abl$ of $7 \times L$ where $L$ is the number of columns of $\Xbl$. For comparison purpose, the number of rows of $\Abl$ is not changed.   All the simulations are conducted in using MATLAB R2016 on Intel(R) Core(TM) i7-6700 CPU@ 3.40 GHz processor. Each optimization problem is repeated five times and the average values are reported.  

Tables~\ref{my-label2} and \ref{my-label1} reports the results of numerical experiments with  the two domains $\mathcal{S}=\{-1,0,1\}$ and $\mathcal{S}=\{-2,-1,0,1, 2\}$.  In Tables~\ref{my-label2} and \ref{my-label1}, $n$ indicates the number of variables to be optimized, and $N_{\mathcal{F},avg}$  the number of average nodes the algorithm visited during the computing set $\mathcal{F}$. For each simulation, the exact solution for $\Xbl$ is obtained. 

The results in Tables~\ref{my-label2}-\ref{my-label1} show that the $N_{\mathcal{F},avg}$  does not depend on the rank of $\Xbl$ but the number of rows of $\Xbl$.  For example, for the number of row equals to ten and the different ranks of  $\Xbl$, $N_{\mathcal{F},avg}$ is same. However, the the number of row increases, the $N_{\mathcal{F},avg}$ also increases. It can be seen that the average CPU time depends on the dimension of $\Xbl$ for the same number of decision variables. For example, $n=100, 300$  the average CPU times increases  with the dimension of $\Xbl$. This is due to the number of  nodes need to be visited  to compute set $\mathcal{F}$. Further, the domains have an effect on the average CPU time dramatically. It can be seen that the average CPU time increased fifty folds for $n=300$  when the domain has extended two more integers.

\begin{table}[]
\centering
\caption{Numerical Results for Problem \eqref{Eq:INMFProblemH} with $\mathcal{S}=\{-1,0,1\}$}
\label{my-label2}
\begin{tabular}{|c|c|c|c|c|}
\cline{1-5}
\begin{tabular}[c]{@{}l@{}}Size of  \\ $\Xbl$\end{tabular} & \begin{tabular}[c]{@{}l@{}}Rank of\\ $\Xbl$\end{tabular} & $n$ & 
 \begin{tabular}[c]{@{}l@{}} Avg CPU \\  time [sec]\end{tabular} & \begin{tabular}[c]{@{}l@{}} $N_{\mathcal{F},avg}$ \end{tabular}  \\\cline{1-5}
                                  10 $\times$ 2 &     2      &    20 &                               0.0859                      & 1767 \\ \cline{1-5}
                                       10 $\times$ 4 &     4                                                     &    40 &                               0.1211                     & 1767 \\ \cline{1-5}
   10 $\times$ 5 &     5                                                     &    50 &                               0.8234                     & 1767 \\ \cline{1-5}
            15 $\times$ 6 &   6 &55&    0.9156  &  19590                             \\ \cline{1-5}
            20 $\times$ 5 &    5 &100&     8.0062 &        411639                          \\ \cline{1-5}
            25 $\times$ 4 &    4 &100&  22.406    &     1534308                       \\ \cline{1-5}
            20 $\times$ 7 &    7 &140&    9.3438   &    411639                                \\ \cline{1-5}
20 $\times$ 8 &   8 &160&     10.5188        &      411639        \\ \cline{1-5}
20 $\times$ 9 &    9 &180& 10.791     &      411639    \\ \cline{1-5}
25 $\times$ 10 &    10 &250&    33.9656
     &     1534308              \\ \cline{1-5}
25 $\times$ 11 &    11 &275&  31.95   &             1534308                       \\ \cline{1-5}
25 $\times$ 12 &   12 &300&    33.1844   &         1534308                          \\ \cline{1-5}
30 $\times$ 10 &    10 &300&     94.3781
  &         7225863                            \\ \cline{1-5}
30 $\times$ 15 &    15 &450&   108.5531    &        7225863                   \\ \cline{1-5}
35 $\times$ 16 &   16 &560&  666.68    &          45800736        \\ \cline{1-5}
40 $\times$ 20 &   20 &800&  2188.5    &          181100487                          \\ \cline{1-5}
45 $\times$ 20 &    20 &900&   9380.1   &        666108006                              \\ \cline{1-5}
50 $\times$ 20 &   20 &1000&   31394.1   &          1703288769                    \\ \cline{1-5}
\end{tabular}
\end{table}

\begin{table}[]
\centering
\caption{Numerical Results for Problem \eqref{Eq:INMFProblemH} with $\mathcal{S}=\{-2,-1,0,1,2\}$}
\label{my-label1}
\begin{tabular}{|c|c|c|c|c|}
\cline{1-5}
\begin{tabular}[c]{@{}l@{}}Size of  \\ $\Xbl$\end{tabular} & \begin{tabular}[c]{@{}l@{}}Rank of\\ $\Xbl$\end{tabular} & $n$ & 
 \begin{tabular}[c]{@{}l@{}}CPU time\\ s\end{tabular} & \begin{tabular}[c]{@{}l@{}} $N_{\mathcal{F},avg}$\end{tabular}  \\\cline{1-5}
                                  10 $\times$ 2 &     2      &    20 &                               0.5828                      & 31815 \\ \cline{1-5}
                                       10 $\times$ 3 &     3                                                     &    30 &                               0.6937                      & 31815 \\ \cline{1-5}
   10 $\times$ 4 &     4                                                     &    40 &                               0.8234                     & 31815 \\ \cline{1-5}
            11 $\times$ 5 &    5 &55&    0.9859  &  31940                             \\ \cline{1-5}
            14 $\times$ 5 &    5 &70&     4.8609 &        240815                          \\ \cline{1-5}
            15 $\times$ 6 &    6 &90&  11.484    &     540400                       \\ \cline{1-5}
            15 $\times$ 7 &    7 &105&   12.083   &    540400                                \\ \cline{1-5}
20 $\times$ 7 &    7 &140&                306.79            &    26713445 \\ \cline{1-5}
20 $\times$ 8 &   8 &160&     318.25        &      26713445        \\ \cline{1-5}
20 $\times$ 10 &    10 &200& 366.36     &      26713445    \\ \cline{1-5}
22 $\times$ 10 &    10 &220& 423.53     &     33838935              \\ \cline{1-5}
22 $\times$ 11 &    11 &242&   488.14    &             33838935                       \\ \cline{1-5}
25 $\times$ 12 &   12 &300&   1244.7   &      97556730                             \\ \cline{1-5}
30 $\times$ 10 &    10 &300&  5096.8    &         466803555                            \\ \cline{1-5}
\end{tabular}
\end{table}

\section{Conclusions}
In this work, we have proposed an exact algorithm to solve a rank-constrained integer least-squares (ILS) problem arising low-rank matrix factorization related applications. The algorithm can also handle the equality constraints and the sparsity constraints. The algorithm is based on the modified Fincke-Pohst enumeration and sparse solutions of Diophantine equations. The proposed algorithm has been explained using an example. Further, numerical experiments have been carried out to  test the applicability of the algorithm.
The results of numerical experiments   show that the proposed algorithm has solved moderate scale rank-constrained integer least-squares problems with the sparsity and equality constraints which arise in intended practical applications.  It is shown that the proposed solution  can be useful to solve  ILS problems with linear equality and sparsity constraints. It has been also observed that the computing sparse solution of Diophantine equations to find row solution sets takes majority time in the algorithm. In future, the proposed algorithm can be improved by reducing  time taken in this step. 
\bibliographystyle{siamplain}
\bibliography{references}

\end{document}

%% file: reverse_sol_tree2.pstex_t
\begin{picture}(0,0)%
\includegraphics{reverse_sol_tree2.pstex}%
\end{picture}%
\setlength{\unitlength}{4144sp}%
\begingroup\makeatletter\ifx\SetFigFont\undefined%
\gdef\SetFigFont#1#2#3#4#5{%
  \reset@font\fontsize{#1}{#2pt}%
  \fontfamily{#3}\fontseries{#4}\fontshape{#5}%
  \selectfont}%
\fi\endgroup%
\begin{picture}(8757,3575)(841,-2833)
\put(856,-1141){\makebox(0,0)[lb]{\smash{{\SetFigFont{20}{24.0}{\rmdefault}{\mddefault}{\updefault}{\color[rgb]{0,0,0}[$x_7$ $x_6$]}%
}}}}
\put(6751,-1186){\makebox(0,0)[lb]{\smash{{\SetFigFont{20}{24.0}{\rmdefault}{\mddefault}{\updefault}{\color[rgb]{0,0,0}[$0$  $0$]}%
}}}}
\put(3736,-1141){\makebox(0,0)[lb]{\smash{{\SetFigFont{20}{24.0}{\rmdefault}{\mddefault}{\updefault}{\color[rgb]{0,0,0}[$-1$  $-1$]}%
}}}}
\put(9226,-1141){\makebox(0,0)[lb]{\smash{{\SetFigFont{20}{24.0}{\rmdefault}{\mddefault}{\updefault}{\color[rgb]{0,0,0}[$1$  $1$]}%
}}}}
\end{picture}%

%% file: reverse_sol_tree4.pstex_t
\begin{picture}(0,0)%
\includegraphics{reverse_sol_tree4.pstex}%
\end{picture}%
\setlength{\unitlength}{4144sp}%
\begingroup\makeatletter\ifx\SetFigFont\undefined%
\gdef\SetFigFont#1#2#3#4#5{%
  \reset@font\fontsize{#1}{#2pt}%
  \fontfamily{#3}\fontseries{#4}\fontshape{#5}%
  \selectfont}%
\fi\endgroup%
\begin{picture}(12907,7321)(66,-6740)
\put(6751,-4291){\makebox(0,0)[lb]{\smash{{\SetFigFont{20}{24.0}{\rmdefault}{\mddefault}{\updefault}{\color[rgb]{0,0,0}[$-1$ $0$ $-1$ $0$ $1$ $1$]}%
}}}}
\put(6661,-196){\makebox(0,0)[lb]{\smash{{\SetFigFont{20}{24.0}{\rmdefault}{\mddefault}{\updefault}{\color[rgb]{0,0,0}[$1$ $0$ $1$ $0$ $-1$ $-1$]}%
}}}}
\put(9991,-196){\makebox(0,0)[lb]{\smash{{\SetFigFont{20}{24.0}{\rmdefault}{\mddefault}{\updefault}{\color[rgb]{0,0,0}[$0$ $1$ $0$ $1$ $0$ $-1$ $-1$]}%
}}}}
\put(6661,-1546){\makebox(0,0)[lb]{\smash{{\SetFigFont{20}{24.0}{\rmdefault}{\mddefault}{\updefault}{\color[rgb]{0,0,0}[$1$ $1$ $-1$ $-1$ $0$ $0$]}%
}}}}
\put(6661,-2311){\makebox(0,0)[lb]{\smash{{\SetFigFont{20}{24.0}{\rmdefault}{\mddefault}{\updefault}{\color[rgb]{0,0,0}[$0$ $0$ $0$ $0$ $0$ $0$]}%
}}}}
\put(6661,-2671){\makebox(0,0)[lb]{\smash{{\SetFigFont{20}{24.0}{\rmdefault}{\mddefault}{\updefault}{\color[rgb]{0,0,0}[$-1$ $1$ $1$ $0$ $0$ $0$]}%
}}}}
\put(6661,-3031){\makebox(0,0)[lb]{\smash{{\SetFigFont{20}{24.0}{\rmdefault}{\mddefault}{\updefault}{\color[rgb]{0,0,0}[$-1$ $-1$ $1$ $1$ $0$ $0$]}%
}}}}
\put(9991,-1501){\makebox(0,0)[lb]{\smash{{\SetFigFont{20}{24.0}{\rmdefault}{\mddefault}{\updefault}{\color[rgb]{0,0,0}[$0$ $1$ $1$ $-1$ $-1$ $0$ $0$]}%
}}}}
\put(6661,-1906){\makebox(0,0)[lb]{\smash{{\SetFigFont{20}{24.0}{\rmdefault}{\mddefault}{\updefault}{\color[rgb]{0,0,0}[$1$ $-1$ $-1$ $0$ $0$ $0$]}%
}}}}
\put(9991,-2266){\makebox(0,0)[lb]{\smash{{\SetFigFont{20}{24.0}{\rmdefault}{\mddefault}{\updefault}{\color[rgb]{0,0,0}[$0$ $0$ $0$ $0$ $0$ $0$ $0$]}%
}}}}
\put(9991,-1861){\makebox(0,0)[lb]{\smash{{\SetFigFont{20}{24.0}{\rmdefault}{\mddefault}{\updefault}{\color[rgb]{0,0,0}[$1$ $1$ $-1$ $-1$ $0$ $0$ $0$]}%
}}}}
\put(9991,-2671){\makebox(0,0)[lb]{\smash{{\SetFigFont{20}{24.0}{\rmdefault}{\mddefault}{\updefault}{\color[rgb]{0,0,0}[$-1$ $-1$ $1$ $1$ $0$ $0$ $0$]}%
}}}}
\put(9991,-3031){\makebox(0,0)[lb]{\smash{{\SetFigFont{20}{24.0}{\rmdefault}{\mddefault}{\updefault}{\color[rgb]{0,0,0}[$0$ $-1$ $-1$ $1$ $1$ $0$ $0$]}%
}}}}
\put(10036,-4246){\makebox(0,0)[lb]{\smash{{\SetFigFont{20}{24.0}{\rmdefault}{\mddefault}{\updefault}{\color[rgb]{0,0,0}[$0$ $-1$ $0$ $-1$ $0$ $1$ $1$]}%
}}}}
\put(1621,-6676){\makebox(0,0)[lb]{\smash{{\SetFigFont{20}{24.0}{\rmdefault}{\mddefault}{\updefault}{\color[rgb]{0,0,0}[$x_6$ $x_7$]}%
}}}}
\put(3286,-6676){\makebox(0,0)[lb]{\smash{{\SetFigFont{20}{24.0}{\rmdefault}{\mddefault}{\updefault}{\color[rgb]{0,0,0}[$x_3$ $x_4$ $x_5$ $x_6$ $x_7$]}%
}}}}
\put(6661,-6676){\makebox(0,0)[lb]{\smash{{\SetFigFont{20}{24.0}{\rmdefault}{\mddefault}{\updefault}{\color[rgb]{0,0,0}[$x_2$ $x_3$ $x_4$ $x_5$ $x_6$ $x_7$]}%
}}}}
\put(9991,-6676){\makebox(0,0)[lb]{\smash{{\SetFigFont{20}{24.0}{\rmdefault}{\mddefault}{\updefault}{\color[rgb]{0,0,0}[$x_1$ $x_2$ $x_3$ $x_4$ $x_5$ $x_6$ $x_7$]}%
}}}}
\put(3331,-5326){\makebox(0,0)[lb]{\smash{{\SetFigFont{20}{24.0}{\rmdefault}{\mddefault}{\updefault}{\color[rgb]{0,0,0}[$0$ $1$ $1$ $1$ $1$]}%
}}}}
\put(3331,-4876){\makebox(0,0)[lb]{\smash{{\SetFigFont{20}{24.0}{\rmdefault}{\mddefault}{\updefault}{\color[rgb]{0,0,0}[$-1$ $0$ $1$ $1$ $1$]}%
}}}}
\put(3331,-4426){\makebox(0,0)[lb]{\smash{{\SetFigFont{20}{24.0}{\rmdefault}{\mddefault}{\updefault}{\color[rgb]{0,0,0}[$1$ $0$ $0$ $1$ $1$]}%
}}}}
\put(3331,-3976){\makebox(0,0)[lb]{\smash{{\SetFigFont{20}{24.0}{\rmdefault}{\mddefault}{\updefault}{\color[rgb]{0,0,0}[$0$ $-1$ $0$ $1$ $1$]}%
}}}}
\put(1666,-4606){\makebox(0,0)[lb]{\smash{{\SetFigFont{20}{24.0}{\rmdefault}{\mddefault}{\updefault}{\color[rgb]{0,0,0}[$1$   $1$]}%
}}}}
\put(1621,-2491){\makebox(0,0)[lb]{\smash{{\SetFigFont{20}{24.0}{\rmdefault}{\mddefault}{\updefault}{\color[rgb]{0,0,0}[$0$   $0$]}%
}}}}
\put(3331,-3256){\makebox(0,0)[lb]{\smash{{\SetFigFont{20}{24.0}{\rmdefault}{\mddefault}{\updefault}{\color[rgb]{0,0,0}[$-1$ $1$ $1$ $0$ $0$]}%
}}}}
\put(3331,-2851){\makebox(0,0)[lb]{\smash{{\SetFigFont{20}{24.0}{\rmdefault}{\mddefault}{\updefault}{\color[rgb]{0,0,0}[$1$ $1$ $0$ $0$ $0$]}%
}}}}
\put(3331,-1996){\makebox(0,0)[lb]{\smash{{\SetFigFont{20}{24.0}{\rmdefault}{\mddefault}{\updefault}{\color[rgb]{0,0,0}[$-1$ $-1$ $0$ $0$ $0$]}%
}}}}
\put(3331,-2401){\makebox(0,0)[lb]{\smash{{\SetFigFont{20}{24.0}{\rmdefault}{\mddefault}{\updefault}{\color[rgb]{0,0,0}[$0$ $0$ $0$ $0$ $0$]}%
}}}}
\put(3331,-1546){\makebox(0,0)[lb]{\smash{{\SetFigFont{20}{24.0}{\rmdefault}{\mddefault}{\updefault}{\color[rgb]{0,0,0}[$1$ $-1$ $-1$ $0$ $0$]}%
}}}}
\put(3286,434){\makebox(0,0)[lb]{\smash{{\SetFigFont{20}{24.0}{\rmdefault}{\mddefault}{\updefault}{\color[rgb]{0,0,0}[$0$ $-1$ $-1$ $-1$ $-1$]}%
}}}}
\put(3286,-61){\makebox(0,0)[lb]{\smash{{\SetFigFont{20}{24.0}{\rmdefault}{\mddefault}{\updefault}{\color[rgb]{0,0,0}[$1$ $0$ $-1$ $-1$ $-1$]}%
}}}}
\put(3286,-466){\makebox(0,0)[lb]{\smash{{\SetFigFont{20}{24.0}{\rmdefault}{\mddefault}{\updefault}{\color[rgb]{0,0,0}[$-1$ $0$ $0$ $-1$ $-1$]}%
}}}}
\put(3286,-916){\makebox(0,0)[lb]{\smash{{\SetFigFont{20}{24.0}{\rmdefault}{\mddefault}{\updefault}{\color[rgb]{0,0,0}[$0$ $1$ $0$ $-1$ $-1$]}%
}}}}
\put(1576,-376){\makebox(0,0)[lb]{\smash{{\SetFigFont{20}{24.0}{\rmdefault}{\mddefault}{\updefault}{\color[rgb]{0,0,0}[$-1$ $-1$]}%
}}}}
\put(10306,-4741){\makebox(0,0)[lb]{\smash{{\SetFigFont{20}{24.0}{\rmdefault}{\mddefault}{\updefault}{\color[rgb]{0,0,0}Final Solution $\mathcal{F}$}%
}}}}
\end{picture}%

%% file: sol_tree1.pstex_t
\begin{picture}(0,0)%
\includegraphics{sol_tree1.pstex}%
\end{picture}%
\setlength{\unitlength}{4144sp}%
\begingroup\makeatletter\ifx\SetFigFont\undefined%
\gdef\SetFigFont#1#2#3#4#5{%
  \reset@font\fontsize{#1}{#2pt}%
  \fontfamily{#3}\fontseries{#4}\fontshape{#5}%
  \selectfont}%
\fi\endgroup%
\begin{picture}(8720,5161)(1831,-4355)
\put(5086,659){\makebox(0,0)[lb]{\smash{{\SetFigFont{12}{14.4}{\rmdefault}{\mddefault}{\updefault}{\color[rgb]{0,0,0}[$x_1$ $x_2$ $x_3$ $x_4$ $x_5$ $x_6$ $x_7$]}%
}}}}
\put(1846,-4291){\makebox(0,0)[lb]{\smash{{\SetFigFont{20}{24.0}{\rmdefault}{\mddefault}{\updefault}{\color[rgb]{0,0,0}$\mathcal{F}_{1,0}$ }%
}}}}
\put(5086,-4246){\makebox(0,0)[lb]{\smash{{\SetFigFont{20}{24.0}{\rmdefault}{\mddefault}{\updefault}{\color[rgb]{0,0,0}$\mathcal{F}_{2,0}$ }%
}}}}
\put(8461,-4201){\makebox(0,0)[lb]{\smash{{\SetFigFont{20}{24.0}{\rmdefault}{\mddefault}{\updefault}{\color[rgb]{0,0,0}$\mathcal{F}_{3,0}$ }%
}}}}
\end{picture}%

%% file: sol_tree2.pstex_t
\begin{picture}(0,0)%
\includegraphics{sol_tree2.pstex}%
\end{picture}%
\setlength{\unitlength}{4144sp}%
\begingroup\makeatletter\ifx\SetFigFont\undefined%
\gdef\SetFigFont#1#2#3#4#5{%
  \reset@font\fontsize{#1}{#2pt}%
  \fontfamily{#3}\fontseries{#4}\fontshape{#5}%
  \selectfont}%
\fi\endgroup%
\begin{picture}(8807,4936)(1744,-4130)
\put(5086,659){\makebox(0,0)[lb]{\smash{{\SetFigFont{14}{16.8}{\rmdefault}{\mddefault}{\updefault}{\color[rgb]{0,0,0}[$x_1$ $x_2$ $x_3$ $x_4$ $x_5$ $x_6$ $x_7$]}%
}}}}
\put(1846,-4021){\makebox(0,0)[lb]{\smash{{\SetFigFont{20}{24.0}{\rmdefault}{\mddefault}{\updefault}{\color[rgb]{0,0,0}$\mathcal{F}_{1,1}$}%
}}}}
\put(5221,-4066){\makebox(0,0)[lb]{\smash{{\SetFigFont{20}{24.0}{\rmdefault}{\mddefault}{\updefault}{\color[rgb]{0,0,0}$\mathcal{F}_{2,1}$}%
}}}}
\put(8641,-4021){\makebox(0,0)[lb]{\smash{{\SetFigFont{20}{24.0}{\rmdefault}{\mddefault}{\updefault}{\color[rgb]{0,0,0}$\mathcal{F}_{3,1}$}%
}}}}
\end{picture}%

%% file: sol_tree3.pstex_t
\begin{picture}(0,0)%
\includegraphics{sol_tree3.pstex}%
\end{picture}%
\setlength{\unitlength}{4144sp}%
\begingroup\makeatletter\ifx\SetFigFont\undefined%
\gdef\SetFigFont#1#2#3#4#5{%
  \reset@font\fontsize{#1}{#2pt}%
  \fontfamily{#3}\fontseries{#4}\fontshape{#5}%
  \selectfont}%
\fi\endgroup%
\begin{picture}(8622,4081)(1786,-3275)
\put(5086,659){\makebox(0,0)[lb]{\smash{{\SetFigFont{14}{16.8}{\rmdefault}{\mddefault}{\updefault}{\color[rgb]{0,0,0}[$x_1$ $x_2$ $x_3$ $x_4$ $x_5$ $x_6$ $x_7$]}%
}}}}
\put(1801,-3211){\makebox(0,0)[lb]{\smash{{\SetFigFont{20}{24.0}{\rmdefault}{\mddefault}{\updefault}{\color[rgb]{0,0,0}$\mathcal{F}_{1,2}$}%
}}}}
\put(5176,-3211){\makebox(0,0)[lb]{\smash{{\SetFigFont{20}{24.0}{\rmdefault}{\mddefault}{\updefault}{\color[rgb]{0,0,0}$\mathcal{F}_{2,2}$}%
}}}}
\put(8416,-3166){\makebox(0,0)[lb]{\smash{{\SetFigFont{20}{24.0}{\rmdefault}{\mddefault}{\updefault}{\color[rgb]{0,0,0}$\mathcal{F}_{3,2}$}%
}}}}
\end{picture}%

%% file: sol_tree4.pstex_t
\begin{picture}(0,0)%
\includegraphics{sol_tree4.pstex}%
\end{picture}%
\setlength{\unitlength}{4144sp}%
\begingroup\makeatletter\ifx\SetFigFont\undefined%
\gdef\SetFigFont#1#2#3#4#5{%
  \reset@font\fontsize{#1}{#2pt}%
  \fontfamily{#3}\fontseries{#4}\fontshape{#5}%
  \selectfont}%
\fi\endgroup%
\begin{picture}(8405,2956)(1921,-2150)
\put(5086,659){\makebox(0,0)[lb]{\smash{{\SetFigFont{14}{16.8}{\rmdefault}{\mddefault}{\updefault}{\color[rgb]{0,0,0}[$x_1$ $x_2$ $x_3$ $x_4$ $x_5$ $x_6$ $x_7$]}%
}}}}
\put(1936,-2086){\makebox(0,0)[lb]{\smash{{\SetFigFont{20}{24.0}{\rmdefault}{\mddefault}{\updefault}{\color[rgb]{0,0,0}$\mathcal{F}_{1,3}$}%
}}}}
\put(5221,-2086){\makebox(0,0)[lb]{\smash{{\SetFigFont{20}{24.0}{\rmdefault}{\mddefault}{\updefault}{\color[rgb]{0,0,0}$\mathcal{F}_{2,3}$}%
}}}}
\put(8326,-2041){\makebox(0,0)[lb]{\smash{{\SetFigFont{20}{24.0}{\rmdefault}{\mddefault}{\updefault}{\color[rgb]{0,0,0}$\mathcal{F}_{3,3}$}%
}}}}
\end{picture}%